\documentclass[11pt]{amsart}
\usepackage{amsmath}
\usepackage{amssymb}
\usepackage{amsfonts}
\usepackage{graphicx}

\DeclareSymbolFont{AMSb}{U}{msb}{m}{n}
\DeclareMathSymbol{\N}{\mathbin}{AMSb}{"4E}
\DeclareMathSymbol{\Z}{\mathbin}{AMSb}{"5A}
\DeclareMathSymbol{\R}{\mathbin}{AMSb}{"52}
\DeclareMathSymbol{\Q}{\mathbin}{AMSb}{"51}
\DeclareMathSymbol{\I}{\mathbin}{AMSb}{"49}
\DeclareMathSymbol{\C}{\mathbin}{AMSb}{"43}

\theoremstyle{definition}

\theoremstyle{corollary}

\theoremstyle{example}

\theoremstyle{note}

\theoremstyle{notation}

\numberwithin{equation}{section}
\begin{document}
\title[Concerning partition regular matrices]
{Concerning partition regular matrices}
\author{Sourav Kanti Patra and Swapan Kumar Ghosh}
\address{Department of Mathematics, Ramakrishna Mission Vidyamandira,
Belur Math, Howrah-711202, West Bengal, India}
\email{souravkantipatra@gmail.com, swapan12345@yahoo.co.in}

\begin {abstract}
A finite or infinite matrix $A$ is image partition regular
provided that whenever $\mathbb N$ is finitely colored, there must
be some $\vec{x}$ with entries from $\mathbb N$ such that all
entries of $A\vec{x}$ are in some color class. In [6] and [8], the
notion of centrally image partition regular matrices were
introduced to extend the results of finite image partition regular
matrices to infinite image partition regular matrices. It was
shown that centrally image partition regular matrices are closed
under diagonal sum. In the present paper, we show that the
diagonal sum of two matrices, one of which comes from the class of
all Milliken-Taylor matrices and the other from a suitable
subclass of the class of all centrally image partition regular
matrices is also image partition regular. This will produce more
image partition regular matrices. We also study the multiple
structures within one cell of a finite partition of $\mathbb N$
$\cdot$

AMS subjclass [2010] : Primary : 05D10 Secondary : 22A15
\end{abstract}

\maketitle

\section {introduction}
In 1933, R. Rado [9] produced a computable characterization,
called the column condition for the (finite) matrices with
rational entries which are kernel partition regular. Kernel
partition regular matrices are those matrices $A$ with the
property that whenever $\mathbb N$ is finitely colored, there
exists some $\vec{x}$ with monochrome entries such that
$A\vec{x}=\vec{0}$. He also extended the result in his later paper
[10] to cover other subsets of $\mathbb R$ (and even of $\mathbb C$).\\

Though several characterizations of (finite) image partition
regular matrices were known, a reasonable characterization of
image partition regular matrices was introduced by Hindman and
Leader in 1993 [4]. A matrix $A$ is said to be image partition
regular if whenever $\mathbb N$ is finitely colored, there will be
some $\vec{x}$ (with entries from $\mathbb N$) such that the
entries of $A\vec{x}$ are monochrome. Image partition regular
matrices generalize many of the classical theorems of Ramsey
Theory. For example, Schur's Theorem [11] and the van der
Waerden's Theorem [12] are equivalent to say that the matrices
$\left(\begin{matrix} 1 & 0 \\ 0 & 1\\ 1 & 1
\end{matrix}\right)$ and for each $n\in \mathbb N$,
$\left(\begin{matrix} 1 & 0 \\ 1 & 1
\\\vdots & \vdots \\ 1 & n-1 \end{matrix}\right)$ are image
partition regular respectively. Some of the characterizations of
finite image partition regular matrices involve the notion of
central sets. Central sets were introduced by Furstenberg and
defined in terms of notion of topological dynamics. A nice
characterization of central sets in terms of algebraic structure
of $\beta\mathbb N$, the Stone-$\breve{C}$ech compactification of
$\mathbb N$ is given in Definition 1.2 (a). Central sets are very
rich in combinatorial properties. The basic fact about central
sets is given by the Central Sets Theorem, which is due to
Furstenberg [3, Proposition 8.21] for the case $S = \mathbb Z$.\\

\textbf{Theorem 1.1.}(Central Sets Theorem) Let S be a semigroup.
Let $\tau$ be the set of sequences $\{y_n\}_{n=1}^\infty$ in $S$.
Let $C$ be a subset of $S$ which is central and let $F\in\mathcal
P_{f}(\tau)$. Then there exist a sequence $\{a_n\}_{n=1}^\infty$
in $S$ and a sequence $\{H_n\}_{n=1}^\infty$ in $\mathcal
P_{f}(\mathbb N)$ such that for each $n\in\mathbb N$, $maxH_n <
minH_{n+1}$ and for each $L\in\mathcal P_{f}(\mathbb N)$ and each
$f\in F$, $\sum_{n\in L}(a_n + \sum_{t\in H_n} f(t))\in C$.\\

We shall present this characterization below, after
introducing the necessary background information.\\

Let $(S,\cdot)$ be an infinite discrete semigroup. Now the points
of $\beta S$ are taken to be the ultrafilters on $S$, the
principal ultrafilters being identified with the points of $S$.
Given $A\subseteq S$ let us set, $\bar{A}=\{p\in\beta S:A\in p\}$.
Then the set $\{\bar{A}:A\subseteq S\}$ will become a basis for a
topology on $\beta S$. The operation $\cdot$ on $S$ can be
extended to the Stone-$\breve{C}$ech compactification $\beta S$ of
$S$ so that $(\beta S,\cdot)$ is a compact right topological
semigroup (meaning that for any $p\in \beta S$, the function
$\rho_p:\beta S \rightarrow \beta S$ defined by $\rho_p(q)=q\cdot
p$ is continuous) with $S$ contained in its topological center
(meaning that for any $x\in S$, the function $\lambda_x:\beta S
\rightarrow \beta S$ defined by $\lambda_x(q)=x\cdot q$ is
continuous). Given $p,q\in\beta S$ and $A\subseteq S, A\in p\cdot
q$ if and only if $\{x\in
S:x^{-1}A\in q\}\in p$, where $x^{-1}A=\{y\in S:x\cdot y\in A\}$.\\

A nonempty subset $I$ of a semigroup $(T,\cdot)$ is called a left
ideal of $T$ if $T\cdot I\subseteq I$, a right ideal if
$I.T\subseteq I$, and a two-sided ideal (or simply an ideal) if it
is both a left and a right ideal. A minimal left ideal is a left
ideal that does not contain any proper left ideal. Similarly, we
can define minimal right ideal and smallest ideal. Any compact
Hausdorff right topological semigroup $(T,.)$ has the unique
smallest two-sided ideal
$$\begin{array}{ccc}
    K(T) & = & \bigcup\{L:L \text{ is a minimal left ideal of } T\} \\
         & = & \,\,\,\,\,\bigcup\{R:R \text{ is a minimal right ideal of } T\}\\
  \end{array}$$

Given a minimal left ideal $L$ and a minimal right ideal $R$ of
$T$, $L\cap R$ is a group, and in particular $K(T)$ contains an
idempotent. An idempotent that belongs to $K(T)$ is called a
minimal idempotent. Given any subset $J\subseteq T$ we shall
use the notation $E(J)$ to denote the set of all idempotents in $J$.\\

\textbf{Definition 1.2.} Let S be a semigroup and let $C\subseteq
S$.

(a) $C$ is called central in $S$ if there is some idempotent $p\in
K(\beta S)$ such that $C\in p$ (Definition 4.42, [7]).

(b) $C$ is called $central^*$ in $S$ if $C\cap A\neq\emptyset$ for
every cental set $A$ in $S$ (Definition 15.3, [7]).\\

Like kernel partition regular matrices all the finite image
partition regular matrices can be described by a computable
condition called the first entries condition. In the following
theorem (Theorem 2.10, [5]) we see that central sets
characterize all finite image partition regular matrices.\\

\textbf{Theorem 1.3.} Let $u,v\in \mathbb N$ and let  $A$ be a $u
\times v$ matrix with entries from $\mathbb Q$. Then the following
statements are equivalent.\\
(a) $A$ is image partition regular.\\
(b) for every additively central subset $C$ of $\mathbb N$, there
exists $\vec{x}\in\mathbb{N}^v$ such that $A\vec{x}\in C^u$.\\

It is an immediate consequence of Theorem 1.3 that whenever $A$
and $B$ are finite image partition regular matrices, so is
$\left(\begin{matrix} A & 0 \\ 0 & B \end{matrix}\right)$, where
$0$ represents a matrix of appropriate size with zero entries.
However, it is a consequence of Theorem 3.14 of [2] that the
corresponding result is not true for infinite image partition
regular matrices. In [6], it was shown that for an image partition
regular matrix $A$, $I(A)$ (Definition 2.1) is a nonempty compact
subset of $(\beta \mathbb N, +)$. It is also a sub-semigroup of
$(\beta \mathbb N,+)$ if $A$ is a finite image partition regular
matrix.\\

In section 2, we have investigated the multiplicative structure of
$I(A)$ in $(\beta \mathbb N,\cdot)$ for finite and infinite matrix
$A$. Using both additive and multiplicative structure of $\beta
\mathbb N$ we show that the diagonal sum of two matrices, one of
which is a Milliken-Taylor matrix and the other is a subtracted
centrally image partition regular matrix (Definition 2.11) is also
an image partition regular matrix.\\

One knows from Finite Sums Theorem (Corollary 5.10, [7]) that
whenever $r\in \mathbb N$ and $\mathbb{N}=\bigcup_{i=1}^r E_i$
there exist $i\in\{1,2,.....,r\}$ and a sequence
$\{x_n\}_{n=1}^\infty$ in $\mathbb N$ such that
$FS(\{x_n\}_{n=1}^\infty)\subseteq E_i$. It is also known that
whenever $r\in \mathbb N$ and $\mathbb{N}=\bigcup_{i=1}^r E_i$
there exist $j\in\{1,2,.....,r\}$ and a sequence
$\{y_n\}_{n=1}^\infty$ in $\mathbb N$ such that
$FP(\{y_n\}_{n=1}^\infty)\subseteq E_j$. The following
theorem (Corollary 5.22, [7]) will allow us to take $i=j$.\\

\textbf{Theorem 1.4.} Suppose $r\in \mathbb{N}$ and
$\mathbb{N}=\bigcup_{i=1}^r E_i$. Then there exist
$i\in\{1,2,.....,r\}$ and sequences $\{x_n\}_{n=1}^\infty$  and
$\{y_n\}_{n=1}^\infty$ in $\mathbb N$ such that
$FS(\{x_n\}_{n=1}^\infty)\bigcup FP(\{y_n\}_{n=1}^\infty)\subseteq
E_i$.\\

We generalize this result in section 3 in matrix version. We also
use both additive and multiplicative structures of $\beta\mathbb
N$ to show that multiple structures induced from matrices with
rational co-efficients will be in one cell of a finite partition
of $\mathbb N$.

\section{Diagonal sum of matrices}

It is natural that the additive structure of $\beta\mathbb N$ is
useful to study the image partition regular matrices over
$(\mathbb N, +)$. In [2], [5], [6] and [8], mostly the additive
structure of $\beta\mathbb N$ is used to study the image partition
regular matrices over $\beta\mathbb N$. In this section, we shall
see that the multiplicative structure of $\beta\mathbb N$ is also
helpful to study the image partition regular matrices over
$(\mathbb N, +)$.\\

We start with the following definition (Definition 2.4, [6]).\\

\textbf{Definition 2.1.} Let $A$ be a finite or infinite matrix
with entries from $\mathbb{Q}$. Then $I(A)=\{p\in\beta\mathbb{N}$
: for every $P\in p$, there exists $\vec{x}$ with entries from
$\mathbb{N}$ such that all entries of $A\vec{x}$ are in $P\}$.\\

We also recall the following lemma (Lemma 2.5, [6]).\\

\textbf{Lemma 2.2.} Let $A$ be a matrix, finite or infinite with
entries from $\mathbb{Q}$.\\
(a) The set $I(A)$ is compact and $I(A)\neq\emptyset$ if and only
if $A$ is image partition regular.\\
(b) If $A$ is finite image partition regular matrix, then $I(A)$
is a sub-semigroup of $(\beta\mathbb{N},+)$.\\

Let us now investigate the multiplicative structure of $I(A)$. In
the following lemma, we shall see that if $A$ is a image partition
regular matrix then $I(A)$ is a left ideal of $(\beta\mathbb
N,\cdot)$. It is also a two-sided ideal of $(\beta\mathbb
N,\cdot)$
provided $A$ is a finite image partition regular matrix.\\

\textbf{Lemma 2.3.} Let $A$ be a matrix, finite or infinite with
entries from
$\mathbb{Q}$.\\
(a) If $A$ is an image partition regular matrix then $I(A)$ is a
left ideal of $(\beta\mathbb{N},\cdot)$.\\
(b) If $A$ is a finite image partition regular matrix then $I(A)$
is a two-sided ideal of $(\beta\mathbb{N},\cdot)$.

\textbf{Proof.}(a) Let $A$ be a $u\times v$ image partition
regular matrix where $u,v\in\mathbb{N}\cup\{\omega\}$. Let
$p\in\beta\mathbb{N}$ and $q\in I(A)$. Also let $U\in p\cdot q$.
Then $\{x\in\mathbb{N} : x^{-1}U\in q\}\in p$. Choose $z\in
\{x\in\mathbb{N} : x^{-1}U\in q\}$. Then $z^{-1}U\in q$. So there
exists $\vec{x}$ with entries from $\mathbb{N}$ such that $y_j\in
z^{-1}U$ for $0\leq j<u$ where
$\vec{y}=A\vec{x}$,\\
$\vec{y}=\left(\begin{matrix} y_0 \\ y_1 \\ y_2 \\
\vdots\end{matrix}\right)$ and
$\vec{x}=\left(\begin{matrix} x_0 \\ x_1 \\ x_2 \\
\vdots\end{matrix}\right)$; $\vec{y}$ and $\vec{x}$ are $u\times
1$ and $v\times 1$ matrices respectively. Now $y_i\in z^{-1}U$ for
$0\leq i<u$ implies that $zy_i\in U$ for $0\leq i<u$. Let
$\vec{x'}=z\vec{x}$ and $\vec{y'}=z\vec{y}$. Then
$\vec{y'}=A\vec{x'}$. So there exists $\vec{x'}$ with entries from
$\mathbb{N}$ such that all entries of $A\vec{x'}$ are in $U$.
Therefore $p\cdot q\in I(A)$. So $I(A)$ is a left ideal of
$(\beta\mathbb{N},\cdot)$.

(b) Let $A$ be a $u\times v$ matrix, where $u,v\in\mathbb{N}$. By
Lemma 2.3(a), $I(A)$ is a left ideal. We now show that $I(A)$ is a
right ideal of $(\beta\mathbb{N},\cdot)$. Let
$p\in\beta\mathbb{N}$ and $q\in I(A)$. Now let $U\in q\cdot p$.
Then $\{x\in\mathbb{N} : x^{-1}U\in p\}\in q$. So there exists
$\vec{x}$ with entries in $\mathbb{N}$ such that
$y_i\in\{x\in\mathbb{N} : x^{-1}U\in p\}$ for $0\leq i<u$,
where $\vec{y}=A\vec{x}$,\\
$\vec{y}=\left(\begin{matrix} y_0 \\ \vdots \\
y_{u-1}\end{matrix}\right)$ and
$\vec{x}=\left(\begin{matrix} x_0 \\ \vdots \\
x_{v-1}\end{matrix}\right)$; $\vec{y}$ and $\vec{x}$ are $u\times
1$ and $v\times 1$ matrices respectively. Now for $0\leq i<u,
y_i\in \{x\in\mathbb{N} : x^{-1}U\in p\}$. Hence $y_i^{-1}U\in p$
for $0\leq i<u$. This implies $\bigcap_{i=0}^{u-1}y_i^{-1}U\in p$.
So $\bigcap_{i=0}^{u-1}y_i^{-1}U\neq\emptyset$. Let $z\in
\bigcap_{i=0}^{u-1}y_i^{-1}U$. Therefore $z\in y_i^{-1} U$ for all
$i\in\{0,1,2,.....,u-1\}$. Hence $y_iz\in U$ for $0\leq i<u$. Let
$\vec{x'}=\vec{x}z$ and $\vec{y'}=\vec{y}z$. Then
$\vec{y'}=A\vec{x'}$. So there exists $\vec{x'}$ with entries from
$\mathbb{N}$ such that all entries of $A\vec{x'}$ are in $U$. Thus
$q\cdot p\in I(A)$. Therefore $I(A)$ is also a right ideal of
$(\beta\mathbb{N},\cdot)$. Hence $I(A)$ is a two-sided ideal
of $(\beta\mathbb{N},\cdot)$.\\

The following theorem was proved in 2003 using combinatorics ([6],
Lemma 2.3). We now provide an alternative proof of this theorem
using the algebra of $(\beta\mathbb N,\cdot)$.\\

\textbf{Theorem 2.4.} Let $A$ and $B$ be finite and infinite image
partition regular matrices respectively (with rational
co-efficients). Then $\left(\begin{matrix} A & 0 \\ 0 & B
\end{matrix}\right)$ is image partition regular.

\textbf{Proof.} Let $r \in \mathbb N$ be given and $\mathbb{N}=
\bigcup_{i=1} ^rE_i$. Suppose that $A$ be a $u \times v$ matrix
where $u,v \in \mathbb N$. Also let $M=\left(\begin{matrix} A & 0 \\
0 & B \end{matrix}\right)$. Now by Lemma 2.3(b), $I(A)$ is a
two-sided ideal of $(\beta\mathbb N,\cdot)$. So $K(\beta\mathbb
N,\cdot) \subseteq I(A)$. Also by Lemma 2.3(a), $I(B)$ is a left
ideal of $(\beta\mathbb N,\cdot)$. Therefore $K(\beta\mathbb
N,\cdot) \cap I(B)\neq \emptyset$. Hence $I(A)\cap
I(B)\neq\emptyset$. Now choose $p\in I(A)\cap I(B)$. Since
$\mathbb{N}= \bigcup_{i=1} ^rE_i$, there exists $k\in
\{1,2,.....,r\}$ such that $E_k\in p$. Thus by definition of
$I(A)$ and $I(B)$, there exist $\vec{x}\in \mathbb{N}^v$ and
$\vec{y}\in \mathbb{N}^\omega$  such that $A\vec{x}\in E_k^u$ and
$B\vec{y}\in E_k^\omega$. Take $\vec{z}=\left(\begin{matrix}
\vec{x} \\ \vec{y}
\end{matrix}\right)$. Then $M\vec{z}=\left(\begin{matrix} A\vec{x}
\\ B\vec{y} \end{matrix}\right)$. So $M\vec{z}\in E_k^\omega$.
Therefore $M = \left(\begin{matrix} A & 0 \\ 0 & B
\end{matrix}\right)$ is image partition regular.\\

The following theorem is the Corollary 2.6 of [6]. We also provide
an alternative proof of this corollary in the following theorem.\\

\textbf{Theorem 2.5.} Let $\mathcal{F}$ denote the set of all
finite image partition regular matrices with entries from $\mathbb
Q$. If $B$ is an image partition regular matrix then
$I(B)\cap(\bigcap_{A\in\mathcal{F}}I(A))\neq\emptyset$.

\textbf{Proof.} Note that for each $A\in\mathcal{F}$, $I(A)$ is a
two-sided ideal of $(\beta\mathbb N,\cdot)$ and therefore
$K(\beta\mathbb N,\cdot)\subseteq I(A)$. Hence $K(\beta\mathbb
N,\cdot)\subseteq\bigcap_{A\in\mathcal{F}}I(A)$. Also by Lemma
2.3(a), $I(B)$ is a left ideal of $(\beta\mathbb N,\cdot)$. Thus
$I(B)\cap K(\beta\mathbb N,\cdot) \neq\emptyset$. Hence
$I(B)\cap(\bigcap_{A\in\mathcal{F}}I(A))\neq\emptyset$.\\

We introduce the following definition to see that the analogous
statements are also true for kernel partition matrices.\\

\textbf{Definition 2.6.} Let $A$ be a finite or infinite matrix
with entries from $\mathbb{Q}$. Then $J(A)=\{p\in\beta\mathbb{N}$
: for every $P\in p$, there exists $\vec{x}$ with entries from
$P$ and $A\vec{x}=\vec{0}\}$.\\

\textbf{Lemma 2.7.} Let $A$ be a matrix, finite or infinite with
entries from $\mathbb{Q}$.\\
(a) The set $J(A)$ is compact and $J(A)\neq\emptyset$ if and only
if $A$ is kernel partition regular.\\
(b) If $A$ is a finite kernel partition regular matrix then $J(A)$
is a sub-semigroup of $(\beta\mathbb{N},+)$.

\textbf{Proof.} Similar proof as given for Lemma 2.5, [6].\\

\textbf{Lemma 2.8.} Let $A$ be a matrix, finite or infinite with
entries from
$\mathbb{Q}$.\\
(a) If $A$ is a kernel partition regular matrix then $J(A)$ is a
left ideal of $(\beta\mathbb{N},\cdot)$.\\
(b) If $A$ is a finite kernel partition regular matrix then $J(A)$
is a two-sided ideal of $(\beta\mathbb{N},\cdot)$.

\textbf{Proof.} Similar proof as given for Lemma 2.3.\\

We now recall the following definition (Definition 1.6(a) of [8]).\\

\textbf{Definition 2.9.} Let $A$ be a $\omega\times\omega$ matrix
with entries from $\mathbb Q$. The matrix $A$ is centrally image
partition regular if for every central subset $C$ of $(\mathbb N,
+)$ there exists $\vec{x}\in\mathbb{N}^\omega$ such that
$A\vec{x}\in C^\omega$.\\

We also recall the following definitions (Definition 1.7 of [5]
and Definition 3.1 of [8]).\\

\textbf{Definition 2.10.} (a) Let $A$ be a $u\times v$ matrix with
rational entries. Then
$A$ is a first entries matrix if \\
(1) no row of $A$ is $\vec{0}$.\\
(2) the first nonzero entry of each row is positive, and\\
(3) the first nonzero entries of any two rows are equal if they
occur in the same column.\\

If $A$ is a first entries matrix and $d$ is the first nonzero
entry of some row of $A$, then $d$ is called a first entry of $A$.
A first entries matrix is said to be monic whenever all the first
entries of $A$ is 1.\\

(b) Let $A$ be a $\omega\times\omega$ matrix with entries from
$\mathbb{Q}.$ Then $A$ is a segmented image partition regular
matrix if

(1) no row of A is $\vec{0}$.

(2) for each $i\in\omega$, $\{j\in\omega:a_{ij}\neq 0\}$ is
finite.

(3) there is an increasing sequence $\{\alpha_n\}_{n=0}^\infty$ in
$\omega$ such that $\alpha_{0}=0$ and for each $ n\in\omega$,
$\{(a_{i,\alpha_n}, a_{i,\alpha_{n+1}},.....,
a_{i,\alpha_{n+1}-1})$ : $i\in\omega\}-\{\vec{0}\}$ is empty or is
the set of rows of a finite image partition regular matrix.\\

If each of these finite image partition regular matrices is a
first entries matrix, $A$ is called a segmented first entries
matrix. If also the first nonzero entry of each $(a_{i,\alpha_n},
a_{i, \alpha_{n+1}},.....,a_{i,\alpha_{n+1}-1})$ is 1, then $A$
is a monic segmented first entry matrix.\\

Here we shall use both the additive and multiplicative structure
of $\beta\mathbb N$ to show that the diagonal sum of two infinite
image partition regular matrices, one of which comes from the
class of all Milliken-Taylor matrices and the other from the class
of all subtracted centrally image partition regular matrices is
also image partition regular. For this we introduce the following
definition.\\

\textbf{Definition 2.11.} Consider the following conditions for a
$\omega\times\omega$ matrix $A$.

(1) no row of $A$ is $\vec{0}$.

(2) for each $i\in\omega$, $\{j\in\omega:a_{ij}\neq 0\}$ is
finite.

(3) If $\vec{c_0},\vec{c_1},\vec{c_2},.....$ be the columns of
$A$, there exist $n\in\omega$ and $k\in\mathbb N$ such that all
the rows of ($\vec{c_{n}}$ $\vec{c_{n+1}}$
$\vec{c_{n+2}}$.....$\vec{c_{n+k-1}}$) are precisely the rows of a
finite image partition regular matrix and the remaining columns
form a centrally image partition regular matrix.

(4) If $\vec{c_0},\vec{c_1},\vec{c_2},.....$ be the columns of
$A$, there exist $n\in\omega$ and $k\in\mathbb{N}$ such that all
the rows of ($\vec{c_{n}}$ $\vec{c_{n+1}}$
$\vec{c_{n+2}}$.....$\vec{c_{n+k-1}}$) are precisely the rows of a
finite image partition regular matrix and the remaining columns of
$A$ form a segmented image partition regular matrix.\\

We call a matrix $A$ to be `subtracted centrally image partition
regular' if it satisfies (1) (2) and (3) and also we call $A$ to
be `subtracted segmented image partition regular' if it satisfies
(1) (2) and (4).\\

Note that both of these matrices defined in Definition 2.11 are
centrally image partition regular. Also note that finite sum
matrices are centrally image partition regular but not subtracted
centrally image partition regular.\\

\textbf{Example 2.12.} An example of a subtracted image partition
regular matrix is given below.\\

$\left(\begin{matrix} 2 & 0 & 1 & 0 & 0 & 0 \cdots \\
                                       2 & 1 & 0 & 0 & 0 & 0 \cdots \\
                                       2 & 1 & 1 & 0 & 0 & 0 \cdots \\
                                       2 & 0 & 0 & 1 & 0 & 0 \cdots \\
                                       2 & 0 & 1 & 0 & 0 & 0 \cdots \\
                                       2 & 1 & 0 & 1 & 0 & 0 \cdots \\
                                       2 & 1 & 1 & 1 & 0 & 0 \cdots \\
                                       2 & 0 & 1 & 1 & 0 & 0 \cdots \\
                                       2 & 0 & 0 & 0 & 0 & 0 \cdots \\
                                       \vdots & \vdots & \vdots & \vdots & \vdots & \vdots
                                       \end{matrix}\right)$\\

We now recall the following definitions (Definition 2.1 of [8],
Definition 5.13(b) of [7] and Definition 2.4 of [2]).\\

\textbf{Definition 2.13.} Let $\vec{x}\in\omega^v$ where
$v\in\mathbb{N}\cup\{\mathbb N\}$.
Then\\
(a) $d(\vec{x})$ is the sequence obtained by deleting all
occurence of $0$ from $\vec{x}$.\\
(b) $c(\vec{x})$ is the sequence obtained by deleting every digit
in $d(\vec{x})$ which is equal to its predecessors and\\
(c) $\vec{x}$ is a compressed sequence if and only if
$\vec{x}=c(\vec{x})$.\\

\textbf{Definition 2.14.} Let $\{x_t\}_{t=0}^\infty$ be a sequence
in $\mathbb{N}$. A sequence $\{y_t\}_{t=0}^\infty$ in $\mathbb{N}$
is said to be a sum-subsystem of $\{x_t\}_{t=0}^\infty$ if there
exists a sequence $\{H_t\}_{t=0}^\infty$ of finite subsets of
$\mathbb{N}$ with $\max H_t<\min H_{t+1}$ for all $t\in\omega$
such that $y_t=\sum_{s\in H_t}x_s$.\\

\textbf{Definition 2.15.} A $\omega\times\omega$ matrix $M$ with
entries from $\omega$ is said to be a Milliken-Taylor matrix if\\
(a) each row of $M$ has only finitely many nonzero entries and\\
(b) there exists a finite compressed sequence $\vec{a}\in\omega^m,
m\in\mathbb{N}$ such that $M$ consists of precisely all the row
vectors $\vec{r}\in\omega^{\omega}$ for which $c(\vec{r})=\vec{a}$.\\

\textbf{Theorem 2.16.} Let  $A$ be a subtracted centrally image
partition regular matrix and $B$ be a Milliken-Taylor matrix. Then
$\left(\begin{matrix} A & 0 \\ 0 & B \end{matrix}\right)$ is image
partition regular.

\textbf{Proof.} If $\vec{c_0},\vec{c_1},\vec{c_2},.....$ be the
columns of $A$ then choose $n,k\in\omega$ and take
$A_1$=($\vec{c_{n}}$ $\vec{c_{n+1}}$.....$\vec{c_{n+k-1}}$) and
$A_2$ to be the remaining columns of $A$ as in the definition of
subtracted centrally image partition regular matrix. Take
$C$=($A_1$ $A_2$) and observe that $I(A)=I(C)$. Now all the rows
of $A_1$ is precisely the rows of a finite image partition regular
matrix (i.e $A_1$ has infinitely many repeated rows precisely
coming from a particular finite image partition regular matrix.).
Therefore $I(A_1)$ is a sub-semigroup of $(\beta\mathbb N,+)$ and
is a two-sided ideal of $(\beta\mathbb N,\cdot)$ by Lemma 2.2(b)
and Lemma 2.3(b) respectively. Also observe that $E(K(\beta\mathbb
N,+)))\subseteq I(A_1)$. Since $B$ is a Milliken-Taylor matrix,
there is a finite compressed sequence $\vec{a} =\{a_i\}_{i=1}^m$,
$m\in\mathbb N$ such that $c(\vec{r})=\vec{a}$ for each row
$\vec{r}$ of $M$. Now choose $p\in E(K(\beta\mathbb N,+))$. Take
$q = a_1\cdot p+a_2\cdot p+.....+a_{m-1}\cdot p$ and $r=a_m\cdot
p$. Then $q\in I(A_1)$ because $p\in E(K(\beta\mathbb N,+))$,
$E(K(\beta\mathbb N,+))\subseteq I(A_1)$ and  $q = a_1\cdot
p+a_2\cdot p+.....+a_{m-1}\cdot p$. Also by Theorem 16.24 of [7],
$cl E(K(\beta\mathbb N,+))$(=$\{p\in\beta\mathbb{N}$ : for each
$P\in p$, $P$ is central in $(\mathbb N,+)\}$) is a left ideal
ideal of $(\beta\mathbb N,\cdot)$. Therefore $r\in cl
E(K(\beta\mathbb N,+))$ as $p\in E(K(\beta\mathbb N,+))\subseteq
cl E(K(\beta\mathbb N,+))$ and $r=a_m\cdot p$. Let $V\in q+r$.
Then $\{x\in\mathbb{N} : -x+V\in r\}\in q$. Now since $q\in
I(A_1)$, there exists $\vec{x^{(1)}}\in \mathbb{N}^k$ such that
$y_i\in\{x\in\mathbb{N} : -x+V\in r\}$ for all $i\in \omega$ where
$\vec{y}=A_1\vec{x^{(1)}}$ and
$\vec{y}=\left(\begin{matrix} y_0 \\ y_1 \\ y_2 \\
\vdots\end{matrix}\right)$. Hence $-y_i+V\in r$ for all
$i\in\omega$. Also observe that $\{y_i : i\in \omega\}$ is finite.
Thus $\bigcap_{i\in\omega}(-y_i+V)\in r$. Now
$\bigcap_{i\in\omega}(-y_i+V)$ is a central set in ($\mathbb N,+$)
because $r\in cl E(K(\beta\mathbb N,+))$ and
$\bigcap_{i\in\omega}(-y_i+V)\in r$. Since $A_2$ is a centrally
image partition regular matrix there exists $\vec{x^{(2)}}\in
\mathbb{N}^{\omega}$ such that whenever $\vec{z}=A_2\vec{x^{(2)}}$
we have $z_j\in\bigcap_{i\in\omega}(-y_i+V)$ for all $j\in\omega$
where $\vec{z}=\left(\begin{matrix} z_0 \\ z_1 \\ z_2 \\
\vdots\end{matrix}\right)$. So $y_i+z_j\in V$ for all
$i,j\in\omega$. Now let $\vec{x}= \left(\begin{matrix}
\vec{x^{(1)}}\\ \vec{x^{(2)}} \end{matrix}\right)$. Then
$C\vec{x}=A_1\vec{x^{(1)}}+A_2\vec{x^{(2)}}=\vec{y}+\vec{z}$.
Therefore $C\vec{x}\in V^\omega$ and hence $q+r\in I(C)=I(A)$.
Also since $p$ is an additive idempotent of $(\beta\mathbb N,+)$,
by Theorem 5.8 of [7] we have a sequence $\{y_t\}_{t=0}^\infty$ in
$\mathbb N$ such that
$p\in\bigcap_{k=1}^\infty\overline{FS(\{y_t\})_{t=k}^{\infty}}$.
Now $V\in a_1\cdot p+a_2\cdot p+.....+a_m\cdot p$ and so by
Theorem 17.31 of [7], there is a sum-subsystem
$\{x_t\}_{t=0}^\infty$ of $\{y_t\}_{t=0}^\infty$ such that
$MT(\vec{a}, \{x_t\}_{t=0}^\infty)\subseteq V$. Let
$\vec{x}=\left(\begin{matrix} x_0 \\ x_1 \\ x_2 \\
\vdots\end{matrix}\right)$. Then $B\vec{x}\in V^\omega$. Thus
$q+r\in I(B)$. Therefore $q+r\in I(A)\bigcap I(B)$. Since
$I(A)\bigcap I(B)\neq\emptyset$, it follows immediately from the
definition of $I(A)$ and $I(B)$ that $\left(\begin{matrix} A & 0 \\
0 & B \end{matrix}\right)$ is image partition regular.\\

Theorem 2.16 is also true if $A$ is taken to be a subtracted
segmented image partition regular matrix.\\

\section{Combined additive and multiplicative structures induced
by matrices}

In this section, we shall present some Ramsey theoretic properties
induced from combined additive and multiplicative structures of
$\mathbb N$. We start with the following theorem (Theorem
15.5 of [7]).\\

\textbf{Theorem 3.1.} Let (S,+) be an infinite commutative
semigroup with identity $0$, let $u,v\in\mathbb{N}$ and let $A$ be
a $u\times v$ matrix with entries from $\omega$ which satisfies
the first entries condition. Let $C$ be a central set in $S$. If
for every first entry $c$ of $A$, $cS$ is a central$^{*}$ set then
there exist sequences $\{x_{1,n}\}_{n=1}^\infty,
\{x_{2,n}\}_{n=1}^\infty,.....,\{x_{v,n}\}_{n=1}^\infty$ in S such
that for every $F\in\mathcal P_f(\mathbb{N}),$\\
$\vec{x}_F\in(S-\{0\})^v$ and $A\vec{x}_F\in C^u$, where
$\vec{x_F}=\left(\begin{matrix} \displaystyle\sum_{n\in F}x_{1,n} \\
\displaystyle\sum_{n\in F}x_{2,n} \\\vdots\\
\displaystyle\sum_{n\in F}x_{v,n}\\
\end{matrix}\right)$.\\

Let $ u,v\in\mathbb{N}\cup\{\omega\}$ and $A$ be a $u \times v$
matrix with entries from $\mathbb{Q}$. Also assume that each row
of $A$ has all but finitely many nonzero entries. Given
$\vec{x}\in \mathbb{N}^v$ and $\vec{y}\in \mathbb{N}^u$, we write
$\vec{x}^A=\vec{y}$ to mean that for $i\in\{0,1,.....,u-1\},$
$\prod_{j=0}^{v-1}x_j^{a_{ij}}=y_i$, where $A= (a_{ij})$ and
$\vec{y}=
\left(\begin{matrix} y_0 \\y_1\\ \vdots \\\end{matrix}\right)$.\\

We now state the following (Theorem 15.20, [7]).\\

\textbf{Theorem 3.2.} Let $u,v\in \mathbb{N}$ and $C$ be a
$u\times v$ matrix with entries from $\mathbb{Z}$. Then the
following statements are equivalent. \\
(a) The matrix $C$ is kernel partition regular over
$(\mathbb{N},+)$. \\
(b) The matrix $C$ is kernel partition regular over
$(\mathbb{N},\cdot)$. That is whenever $r\in \mathbb{N}$ and
$\mathbb{N}-\{1\}=\bigcup_{i=1} ^r D_i$, there exist
$i\in\{1,2,.....,r\}$ and $\vec{x}\in(D_i)^v$ such that
$\vec{x}^C=\vec{1}$.\\
(c) The matrix $C$ satisfies the columns condition over $\mathbb
Q$.\\

Note that the matrix $A=(2,-2,1)$ satisfies the columns condition
over $\mathbb{Q}$. So by Theorem 3.2, $A$ is kernel partition
regular over $(\mathbb N,\cdot)$. Now $D=\{x^2:x\in \mathbb{N}\}$
is not central in $(\mathbb{N},\cdot)$. If $C=\mathbb{N}-D$ then
$C$ is central in $(\mathbb{N},\cdot)$ containing no solution of
the equation $\vec{x}^A=\vec{1}$. We know that image partition
regularity and kernel partition regularity are dual to each other
for finite matrices. So we can find an image partition regular
matrix $B = (2)$ say, containing no image of the matrix $B$ to
some central subset of $(\mathbb{N},\cdot)$. Lemma 3.4 shows which
image partition regular matrix over $(\mathbb N,+)$ contains
solution in central subsets of $(\mathbb{N},\cdot)$.\\

Let $A\subseteq\mathbb{N}$. It is trivial that $A$ is
multiplicatively thick if and only if $A\cap A/2\cap
A/3\cap.....\cap A/n\neq\emptyset$ for each $n\in \mathbb{N}$.
Thick sets are central as follows from Lemma
5.10 of [1], which we state in the following.\\

\textbf{Lemma 3.3.} Any additively thick set in $\mathbb{N}$ is
any additively central and any multiplicatively thick set in
$\mathbb{N}$ is multiplicatively central.\\

\textbf{Lemma 3.4.} For each $k\in \mathbb{N}$ let $A_k=\{x^k:x\in
\mathbb{N}\}$. Then $A_k$ is $central^*$ in $(\mathbb N,\cdot)$ if
and only if $k=1$.

\textbf{Proof.} Let $k\neq1$ and $B_k=\mathbb{N}-A_k$. Choose
$n\in\mathbb{N}$ and select a prime $p > n$. Then $ip\notin A_k$
for $i\in\{1,2,.....,n\}$ and hence $ip\in B_k$ for
$i\in\{1,2,.....,n\}$. Thus $p\in B_k/i$ for
$i\in\{1,2,.....,n\}$. So $p\in B_k\cap B_k/2\cap.....\cap B_k/n$
and thus $B_k\cap B_k/2\cap.....\cap B_k/n\neq\emptyset$.
Consequently, $B_k$ is multiplicatively thick and therefore $B_k$
is multiplicatively central (Lemma 3.3).
Thus $A_k$ is not central$^{*}$ in $(\mathbb N,\cdot)$.\\

If $A$ is a `segmented first entry matrix' with entries from
$\omega$, then one may expect to have a combined additive and
multiplicative structure induced by $A$ in a cell of any finite
coloring of $\mathbb N$. Lemma 3.4 motivates us to consider $A$ to
be a `monic segmented first entry matrix'.\\

\textbf{Lemma 3.5.} Let $A$ be a monic segmented first entry
matrix with entries from $\omega$ and let $C$ be a central subset
of ($\mathbb N,\cdot$). Then there exists
$\vec{x}\in\mathbb{N}^{\omega}$ such that $\vec{x}^A\in
C^{\omega}$.

\textbf{Proof.} Let $\vec{c_{0}}, \vec{c_{1}}, \vec{c_{2}},.....$
denote the the columns of $A$. Let $\{\alpha_n\}_{n=0}^\infty$ be
an increasing sequence as in the definition of a monic segmented
first entry matrix (Definition 2.10 (b)). For each $n\in \omega$,
let $A_{n}$ be the matrix whose columns are $\vec{c_{\alpha_{n}}},
\vec{c_{\alpha_{n}+1}}, \vec{c_{\alpha_{n}+2}},.....,
\vec{c_{\alpha_{n+1}-1}}$. Then the set of nonzero rows of $A_n$
is finite and, if nonempty, is the set of rows of a finite monic
first entries matrix. ($A_n$ may contain infinitely many nonzero
rows but only finitely many rows are distinct.) Let $B$=($A_0$
$A_1$ ..... $A_n$). Let $C$ be a central sunset of ($\mathbb
N,\cdot$) and $p$ be a minimal idempotent in ($\beta\mathbb
N,\cdot$) such that $C\in p$. Let $C^*=\{n\in C:n^{-1}C\in p\}$.
Then $C^*\in p$ and for every $n\in C^*$, $n^{-1}C^*\in p$ (by
Lemma 4.14, [7]). Now by Theorem 3.1, we can choose
$\vec{x^{(0)}}\in \mathbb N^{\alpha_{1}-\alpha_{0}}$ such that if
$\vec{y}=\vec{x^{(0)}}^{A_{0}}$ then $y_i\in C^*$ for every
$i\in\omega$ for which $i^{th}$ row of $A_0$ is nonzero. We now
make induction assumption that, for some $m\in \omega$, we have
chosen $\vec{x^{(0)}}$, $\vec{x^{(1)}}$,.....,$\vec{x^{(m)}}$ such
that if $i\in\{0,1,2,.....,m\}$ then $\vec{x^{(i)}}\in \mathbb
N^{\alpha_{i+1}-\alpha_{i}}$ and if
$\vec{y}=\left(\begin{matrix} \vec{x^{(0)}} \\ \vec{x^{(1)}}  \\
\vdots \\ \vec{x^{(m)}}\end{matrix}\right)^{B_m}$ then $y_j\in
C^*$ for every $j\in \omega$ for which the $j^{th}$ row of $B_m$
is nonzero. Let $D=\{j\in \omega : j^{th}$ row of $B_m$ is
nonzero$\}$ and note that for each $j\in\omega$, $y_j^{-1}C^*\in
p$. By Theorem 3.1, we can choose $\vec{x^{(m+1)}}\in \mathbb
N^{\alpha_{m+2}-\alpha_{m+1}}$ such that if
$\vec{z}=\vec{x^{(m+1)}}^{A_{m+1}}$ then $z_j\in
C^*\cap\bigcap_{t\in D}y_t^{-1}C^*$ for every $j\in \omega$ for
which $j^{th}$ row of $A_{m+1}$ is nonzero and is equal to 1
otherwise. Thus we can choose an infinite sequence
$\{\vec{x^{(i)}}\}_{i\in \omega}$ such that if $i\in \omega$ then
$\vec{x^{(i)}}\in \mathbb N^{\alpha_{i+1}-\alpha_{i}}$
and if $\vec{y}=\left(\begin{matrix} \vec{x^{(0)}} \\
\vec{x^{(1)}}  \\ \vdots \\
\vec{x^{(i)}}\end{matrix}\right)^{B_i}$ then $y_j\in C^*$ for
every $j\in \omega$ for which the $j^{th}$ row of $B_m$ is
nonzero. Let
$\vec{x}=\left(\begin{matrix} \vec{x^{(0)}} \\ \vec{x^{(1)}} \\
\vec{x^{(2)}}  \\ \vdots \end{matrix}\right)$ and let
$\vec{y}=\vec{x}^A$. We note that for every $j\in \omega$, there
exists $m\in \omega$ such that $y_j$ is the $j^{th}$ entry of
$\left(\begin{matrix} \vec{x^{(0)}} \\ \vec{x^{(1)}}  \\ \vdots \\
\vec{x^{(i)}}\end{matrix}\right)^{B_i}$ whenever $i>m$. Thus all
the entries of $\vec{y}$ are in $C^{*}$.\\

\textbf{Theorem 3.6.} Let $A$ be a $\omega\times\omega$ monic
segmented first entries matrix with entries from $\omega$. Let $
r\in\mathbb{N}$ and $\mathbb{N}= \bigcup_{i=1} ^r E_i$ be an
$r$-coloring of $\mathbb{N}$. Then there exist
$i\in\{1,2,.....,r\}$ and vectors $\vec{x},
\vec{y}\in\mathbb{N}^{\omega}$ such that all the elements of
$A\vec{x}$ and $\vec{y}^A$ are in $E_i$.

\textbf{Proof.} Since $I(A)$ is a left ideal of
$(\beta\mathbb{N},\cdot)$, we can choose $p\cdot p = p\in
K(\beta\mathbb{N},\cdot)\bigcap I(A)$. Also since $\mathbb{N}=
\bigcup_{i=1} ^nE_i$ there exists $i\in\{1,2,.....,r\}$ such that
$E_i\in p$. As $E_i\in p$ and $p\in I(A)$, we can choose
$\vec{x}\in \mathbb{N}^{\omega}$ such that $A\vec{x}\in
E_i^{\omega}$. Also $E_i\in p$ and $p\in K(\beta\mathbb{N},\cdot)$
imply that  $E_i$ is a multiplicatively central in ($\mathbb
N,\cdot$). Now by previous lemma, we can find
$\vec{y}\in\mathbb{N}^{\omega}$ such that $\vec{y}^A\in
E_i^{\omega}$. So our claim is proved.\\

We note that the Theorem 1.4 follows as a corollary of the above
theorem.\\

\textbf{Corollary 3.7.} Let $ u,v \in \mathbb{N}$ and $A$ be a
$u\times v$ monic first entry matrix with entries from $\omega$.
Let $r\in \mathbb{N}$ and $\mathbb{N}= \bigcup_{i=1} ^r E_i$ be an
$r$-coloring of $\mathbb{N}$. Then there exist $i \in
\{1,2,.....,r\}$ and $\vec{x},\vec{y}\in \mathbb{N}^v$ such that
$A\vec{x}\in E_i^u$ and $\vec{y}^A\in E_i^u$.

\textbf{Proof.} The proof is almost same as that of Theorem 3.6.\\

We now raise the following question.\\

\textbf{Question 3.8.} Let $ u,v \in \mathbb{N}$ and let $A$ be a
$u\times v$ monic first entry matrix with entries from $\omega$.
Let $r\in \mathbb{N}$ and $\mathbb{N}= \bigcup_{i=1} ^r E_i$ be an
$r$-coloring of $\mathbb{N}$. Does there exist $\vec{x}\in
\mathbb{N}^u$ such that $A\vec{x}\in E_i^u$ and $\vec{x}^A \in
E_i^u$ for some $i \in \{1,2,.....,r\}$?\\

If the question 3.8 is true, one may extend it by taking $u,
v\in\mathbb N\cup\{\omega\}$ by considering $A$ to be a monic
segmented first entry matrix.\\

\textbf{Corollary 3.9.} Let $u,v,r\in\mathbb{N}$ and let $A$ be a
$u\times v$ matrix with entries from $\mathbb{Q}$ satisfying
columns condition over $\mathbb{Z}$. If $\mathbb{N} =
\bigcup_{i=1} ^r E_i$ then there exist $i \in \{1,2,.....,r\}$ and
$\vec{x},\vec{y}\in E_i^u$ such that  $A\vec{x}=\vec{0}$ and
$\vec{y}^A=\vec{1}$.

\textbf{Proof}. Imitate the proof of Theorem 3.6 using Lemma 2.8
and Theorem 15.16(a), [7].\\

We devote the remaining portion of this section to investigate
some Ramsey theoretic properties concerning product of sums or sum
of products arising from matrices. Henceforth, all the matrices
under consideration are with rational entries.\\

Note that $I(A)$ contains all minimal idempotents of
$(\beta\mathbb N,+)$ for all $A\in\mathcal{F}$ where $\mathcal{F}$
denotes the set of all finite image partition regular matrices
over $\mathbb Q$. Hence $I=\bigcap_{A\in\mathcal{F}}I(A)$ contains
all the minimal idempotents of $(\beta\mathbb N,+)$. In Theorem
3.17, we show that multiple structures induced by image partition
regular matrices are contained in one cell of a finite coloring of
$\mathbb N$.\\

\textbf{Definition 3.10.} Let $m\in\mathbb{N}$ be given and
$\vec{y^{(i)}}\in\mathbb{N}^{u_i}$,
$u_i\in\mathbb{N}\cup\{\omega\}$ for $i\in\{1,\cdots ,m\}$. Also
let $\{x_t\}_{t=1}^{\infty}$ be a sequence in $\mathbb{N}$. Then\\
(a) $P(\vec{y^{(1)}},\vec{y^{(2)}},.....,
\vec{y^{(m)}})=\{\prod_{i=1}^{m} y_i:y_i\in\vec{y^{(i)}}$ for all
$i = 1,2,....,m\}$.\\
(b) $S(\vec{y^{(1)}},\vec{y^{(2)}},.....,
\vec{y^{(m)}})=\{\sum_{i=1}^{m} y_i:y_i\in\vec{y^{(i)}}$ for all
$i = 1,2,....,m\}$.\\
(c) $PS_m(\{x_t\}_{t=1}^\infty)=\{\prod_{i=1}^{m}\displaystyle
\sum_{t\in F_{i}}x_t:F_1,F_2,.....,F_m\in
\mathcal{P}_f(\mathbb{N})$
and $maxF_i < minF_{i+1}$ for all $i = 1,....,m-1\}$.\\

\textbf{Lemma 3.11.} Let $m\in\mathbb{N}$ and for each
$i\in\{1,2,.....,m\}$, let $A_i$ be a $u_i\times v_i$ image
partition regular matrix over $\mathbb{N}$ where
$u_i,v_i\in\mathbb{N}$. If $U\in p^m$ where $p\in
I=\bigcap_{A\in\mathcal{F}}I(A)$, there exists
$\vec{x^{(i)}}\in\mathbb{N}^{u_i}$, $1\leq i\leq m$ such that
$P(A_1\vec{x^{(1)}},A_2\vec{x^{(2)}},.....,A_m\vec{x^{(m)}})\subseteq
U$.

\textbf{Proof.} We shall prove this theorem by induction on $m$.
Clearly, by definition of $I$, the theorem is true for $m=1$. Let
the theorem be true for $m=n$. Let $U\in p^{n+1}$ and for each
$i\in\{1,2,.....,n+1\}$, let $A_i$ be a $u_i\times v_i$ image
partition regular matrix over $\mathbb{N}$. Now $U\in p^{n+1}$
implies that $\{x\in\mathbb{N} : x^{-1}U\in p\}\in p^n$. Thus by
induction hypothesis, there exists
$\vec{x^{(i)}}\in\mathbb{N}^{v_i}$ for $1\leq i \leq n$ such that
$P(A_1\vec{x^{(1)}},A_2\vec{x^{(2)}},.....,A_n\vec{x^{(n)}})\subseteq
\{x\in\mathbb{N} : x^{-1}U\in p\}$. Let
$\vec{y^{(i)}}=A_i\vec{x^{(i)}}$ for $1\leq i \leq n$. Then
$(\prod_{i=1}^{n} y_i)^{-1}U\in p$ whenever $y_i\in\vec{y^{(i)}}$
for $1\leq i \leq n$. Let $Y=\{\prod_{i=1}^{n}
y_i:y_i\in\vec{y^{(i)}}$ for $1\leq i \leq n\}$. Since
$\vec{y^{(i)}}$ is finite for each $i\in\{1,2,.....,n\}$, $Y$ is
also finite. Thus we have $\bigcap_{y\in Y}y^{-1}U\in p$. Now
$I\subseteq I(A_{n+1})$. Hence $p\in I(A_{n+1})$. Also since
$\bigcap_{y\in Y}y^{-1}U\in p$, there exists $\vec{x^{(n+1)}}\in
\mathbb{N}^{u_{n+1}}$ such that $y_{n+1}\in\bigcap_{y\in
Y}y^{-1}U$ for all $y_{n+1}\in\vec{y^{(n+1)}}$ where
$\vec{y^{(n+1)}}=A\vec{x^{(n+1)}}$. It is easy to see that
$\prod_{i=1}^{n+1} y_i\in U$ for all $y_i \in\vec{y^{i}}$ for $1
\leq i \leq n+1$. So
$P(\vec{y^{(1)}},\vec{y^{(2)}},.....,\vec{y^{(n+1)}})\subseteq U$.
Therefore
$P(A_1\vec{x^{(1)}},A_2\vec{x^{(2)}},.....,A_{n+1}\vec{x^{(n+1)}})
\subseteq U$. This completes the proof.\\

Note that in above theorem, we need not assume $p$ to be a minimal
idempotent of $(\beta\mathbb{N},+)$. We now prove a similar
version of Lemma 3.11 by replacing one of the finite image
partition regular matrices by an infinite image partition regular
matrix.\\

\textbf{Lemma 3.12.} Let $m\in\mathbb N$ be given and let for each
$i\in\{1,2,.....,m\}$, $A_i$ be a $u_i\times v_i$ image partition
regular matrix where $u_i,v_i\in\mathbb{N}$. Also let $B$ be any
infinite image partition regular matrix. If $U\in p^m\cdot q$
where $p\in I=\bigcap_{A\in\mathcal{F}}I(A)$ and $q\in I(B)$ then
there exist $\vec{x_i}\in \mathbb{N}$ for each
$i\in\{1,2,.....,m\}$ and $\vec{x^{(m+1)}}\in \mathbb{N}^\omega$
such that $P(A_1\vec{x^{(1)}},A_2\vec{x^{(2)}},.....,
A_m\vec{x^{(m)}},B\vec{x^{(m+1)}})\subseteq U$.

\textbf{Proof.} Let $U\in p^m\cdot q$. Then $\{x\in\mathbb{N} :
x^{-1}U\in q\}\in p^m$. By Lemma 3.11, there exists
$\vec{x^{(i)}}\in\mathbb{N}^{u_i}$, $1\leq i\leq m$ such that
$P(A_1\vec{x^{(1)}},A_2\vec{x^{(2)}},.....,A_m\vec{x^{(m)}})\subseteq
\{x\in\mathbb{N} : x^{-1}U\in q\}$. Now let for each
$i\in\{1,2,.....,m\}$, $y_i=A\vec{x^{(i)}}$. For simplicity, let
$Y=P(A_1\vec{x^{(1)}},A_2\vec{x^{(2)}},.....,A_m\vec{x^{(m)}})$.
Since $A_i$ is a finite image partition regular matrix for each
$i\in\{1,2,.....,m\}$, $Y$ is finite. Thus we have $\bigcap_{y\in
Y}y^{-1}U\in q$. Also since $q\in I(B)$, there exists
$\vec{x^{(m+1)}}\in \mathbb{N}^\omega$ such that
$y_{m+1}\in\bigcap_{y\in Y}y^{-1}U$ for all
$y_{m+1}\in\vec{y^{(m+1)}}$ where
$\vec{y^{(m+1)}}=B\vec{x^{(m+1)}}$. Hence $yy_{m+1}\in U$ for all
$y\in Y$ and $y_{m+1}\in\vec{y^{(m+1)}}$. Therefore
$P(A_1\vec{x^{(1)}},A_2\vec{x^{(2)}},.....,A_m\vec{x^{(m)}},
B\vec{x^{(m+1)}}) \subseteq U$.\\

\textbf{Corollary 3.13.} Let $r,m\in\mathbb{N}$ and for each
$i\in\{1,2,.....,m\}$, let $A_i$ be a $u_i\times v_i$ image
partition regular matrix over $\mathbb{N}$;
$u_i,v_i\in\mathbb{N}$. If $\mathbb{N}=\bigcup_{i=1}^{r}E_i$ then
there exist $k\in\{1,2,.....,r\}$ and $\vec{x^{(i)}}\in\mathbb
N^{v_i}$ for $1\leq i\leq m$ such that
$P(A_1\vec{x^{(1)}},A_2\vec{x^{(2)}},.....,
A_m\vec{x^{(m)}})\subseteq E_k$.

\textbf{Proof.} Let $p$ be a minimal idempotent of
$(\beta\mathbb{N},+)$. Choose $k\in\{1,2,.....,r\}$ such that
$E_k\in p^m$ and use Lemma 3.11.\\

Note that it follows by Lemma 3.12 that one can take one of the
matrices $A_1,A_2,.....,A_m$ to be infinite image partition
regular in Corollary 3.13. Theorem 3.14 is an analogous version of
Corollary 3.13 for the sums of products induced by the certain
class of image partition regular matrices.\\

\textbf{Theorem 3.14.} Let $m,r \in \mathbb{N}$ be given and for
each $i\in\{1,2,.....,m\}$, let $A_i$ be a $u_i\times v_i$ monic
first entry matrix, $u_i,v_i\in\mathbb{N}$. If
$\mathbb{N}=\bigcup_{i=1}^{r}E_i$ then there exists $\vec{x^{(i)}}
\in \mathbb{N}^{u_i}$, $1 \leq i \leq m$ such that
$S(\vec{x^{(1)}}^{A_1},\vec{x^{(2)}}^{A_2},.....,
\vec{x^{(m)}}^{A_m})\in E_i$ for some $i\in \{1,2,.....,m\}$.

\textbf{Proof.} See proof of Lemma 3.11 and Corollary 3.13.\\

\textbf{Theorem 3.15.} Let $p+p=p\in\beta\mathbb{N}$. Let
$m\in\mathbb{N}$ and let $U\in p{^m}$. Then there exists a
sequence $\{x_t\}_{t=1}^\infty$ in $\mathbb{N}$ such that
$PS_m(\{x_t\}_{t=1}^\infty)\subseteq U$.

\textbf{Proof.} Imitate the proof of Theorem 17.24 of [7].\\

In the following theorem one does not require $p$ to be an
additive idempotent.\\

\textbf{Theorem 3.16.} Let $\{x_t\}_{t=1}^\infty$ be a sequence in
$\mathbb{N}$. If
$p\in\bigcap_{k=1}^\infty\overline{FS(\{x_t\}_{t=k}^\infty)}$,
then for all $n$ and $k$ in $\mathbb{N}$,
$PS_m(\{x_t\}_{t=1}^\infty) \in p^n$.

\textbf{Proof.} Imitate the proof of Theorem 17.26, [7].\\

\textbf{Theorem 3.17.} Let $r,m\in\mathbb{N}$ be given and
$\mathbb{N}= \bigcup_{i=1} ^nE_i$. Suppose for each
$i\in\{1,2,.....,m\}$, $A_i$ is a $u_i\times v_i$ image partition
regular matrix over $(\mathbb{N},+)$ where $u_i,v_i\in
\mathbb{N}$. Then there exist $k\in\{1,2,.....,m\}$;
$\vec{x^{(i)}}, \vec{y^{(i)}}\in \mathbb{N}^v$ for
$i\in\{1,2,.....,m\}$ and $\{z_t\}_{t=1}^\infty$ such that
$A_i\vec{x^{(i)}}\in E_k^{u_i}$ for $i\in\{1,2,.....,m\}$ and
$P(A_1\vec{y^{(1)}},A_2\vec{y^{(2)}},.....,A_m\vec{y^{(m)}})$,
$PS_m(\{z_t\}_{t=1}^\infty) \subseteq E_k$.

\textbf{Proof.} For each $i\in\{1,2,.....,m\}$, $A_i$ is a finite
image partition regular matrix over $(\mathbb{N},+)$ and hence by
Lemma 2.2, $I(A_i)\neq\emptyset$ and is a compact sub-semigroup of
$(\beta\mathbb{N},+)$ for each $i\in\{1,2,.....,m\}$. Also observe
that $E(K(\beta\mathbb{N},+))\subseteq I(A_i)$ for each
$i\in\{1,2,.....,m\}$. Therefore $\bigcap_{i=1} ^m
I(A_i)\neq\emptyset$ and is a sub-semigroup of
$(\beta\mathbb{N},+)$. Now choose an idempotent $p$ of
$(\beta\mathbb{N},+)$ such that $p\in\bigcap_{i=1} ^m I(A_i)$.
Also $\bigcap_{i=1} ^m I(A_i)$ is a two-sided ideal of
$(\beta\mathbb{N},\cdot)$ because $I(A_i)$ is a two-sided ideal of
$(\beta\mathbb{N},\cdot)$ for each $i\in\{1,.....,m\}$ by Lemma
2.3(b). Thus $p^m\in\bigcap_{i=1} ^m I(A_i)$. As $\mathbb{N}=
\bigcup_{i=1} ^nE_i$, choose $E_k\in p^m$
for some $k\in\{1,2,.....,r\}$. Then\\
(a) For each $i\in\{1,2,.....,m\}$ there exists $\vec{x^{(i)}}\in
\mathbb{N}^v$ such that $A_i\vec{x^{(i)}}\in E_k^{u_i}$, because
$E_k\in p^m$ and $p^m\in\bigcap_{i=1} ^m I(A_i)$.\\
(b) Since $E_k \in p^m$, by Lemma 3.11, for each
$i\in\{1,2,.....,m\}$ there exists $\vec{y_i}\in \mathbb{N}^v_i$
such that
$P(A_1\vec{y^{(1)}},A_2\vec{y^{(2)}},.....,A_m\vec{y^{(m)}})
\subseteq E_k$.\\
(c) As  $E_k \in p^m$, by Theorem 3.15, there exists a sequence
$\{z_t\}_{t=1}^\infty$ such that $PS_m(\{z_t\}_{t=1}^\infty)
\subseteq E_k$. This complete completes the  proof.\\

Similar result is also true if we replace
$P(A_1\vec{y^{(1)}},A_2\vec{y^{(2)}},.....,A_m\vec{y^{(m)}})$ by
$S(\vec{y^{(1)}}^{A_1},\vec{y^{(2)}}^{A_2},.....,
\vec{y^{(m)}}^{A_m})$ and $PS_m(\{z_t\}_{t=1}^\infty)$ by
$SP_m(\{z_t\}_{t=1}^\infty)$ in previous theorem.\\

Note that Lemma 3.12 allows us to take one of the matrices to be
infinite image partition regular over $(\mathbb N, +)$ in the
previous theorem. At the end of our paper, we raise the following
two questions.\\

\textbf{Question 3.18.} Suppose $m\in \mathbb N$ be given and let
for each $i\in \{1,2,.....,m\}$, $A_i$ be a $u_i\times v$ image
partition regular matrix where $u_i,v\in\mathbb{N}$. Let
$r\in\mathbb{N}$ and $\mathbb{N}=\bigcup_{i=1}^{r}E_i$. Does there
exist $k\in\{1,2,.....,r\}$ and $\vec{x}\in\mathbb{N}^v$ such that
$P(A_1\vec{x},A_2\vec{x},.....,A_m\vec{x})\subseteq E_k$?\\

\textbf{Question 3.19.} Let $m\in\mathbb{N}-\{1\}$. Does there
exist a finite partition $\mathcal{R}$ of $\mathbb N$ such that
given any $A\in\mathcal {R}$ there do not exist one-to-one
sequences $\{x_t\}_{t=1}^\infty$ and $\{y_t\}_{t=1}^\infty$ with
$PS_m(\{x_t\}_{t=1}^\infty)\subseteq A$ and
$PS_1(\{y_t\}_{t=1}^\infty)\subseteq A$?\\

\textbf {Acknowledgement.}  The authors are grateful to Prof.
Dibyendu De of Kalyani University for continuous inspiration and a
number of valuable suggestions towards the improvement of the
paper.

\end{document}